\documentclass[10pt]{amsart}
\usepackage{amsthm}
\usepackage[cyr]{aeguill}
\usepackage{dsfont}
\usepackage{upgreek}
\usepackage{mathrsfs}
\usepackage{amsmath,amsfonts,amssymb,color,multicol,subeqnarray}\usepackage{graphicx,longtable}
\usepackage[T1]{fontenc}
\usepackage[latin1]{inputenc}  
\usepackage[a4paper]{geometry}
\theoremstyle{plain}
\newtheorem{teo}{Theorem}
\newtheorem{pro}{Proposition}
\newtheorem{lem}{Lemma}
\newtheorem{cor}{Corollary}
\theoremstyle{definition}
\newtheorem{defi}{Definition}
\theoremstyle{remark}
\newtheorem*{nota}{Remark}
\def\sca #1_#2_#3 {\langle#1,#2\rangle_#3}
\def\R{\mathbb{R}}
\def\res {\textit{\textbf{Res}}}
\def\C{\mathbb {C}}
\def\m{\mu}
\def\<{\langle}
\def\>{\rangle}
\def\bi{\natural}
\def\n{\nu}
\def\pii{\mathcal{\pi}}
\def\hphi{{\mathcal H}^\varphi}
\def\hphik{{\mathcal H}_{K}^\varphi}
\def\pphi{\mathcal{\pi}^\varphi}

\def\Pin{\mathcal{P}_{\leq 1}^\bi(G)}

\def\Pun{\mathcal{P}_1^\bi(G)}

\def\pp{\varphi}
\def\H{{\mathcal H}}
\def\P{{\mathcal {P}}}

\def\A{{\mathcal A}}

\def\E{{\mathcal E}}
\def\Q{{\mathcal Q}}
\def\L{{\mathcal L}}

\def\M{{\mathcal M}}

\newcommand{\bq}{\begin{equation}}

\newcommand{\eq}{\end{equation}}
\newcommand{\ba}{\begin{eqnarray*}}
\newcommand{\ea}{\end{eqnarray*}}

\newcommand{\into}{\int_{G}}

\newcommand{\ove}[1]{\overline{#1}}

\newcommand{\ban}{\begin{eqnarray}}
\newcommand{\ean}{\end{eqnarray}}

\title{A Bochner type theorem for inductive limits of Gelfand pairs} 
\date{October 13, 2007 }
 
\begin{document}
\maketitle
\centerline{M. RABAOUI\footnote{Laboratoire de Math{\'e}matiques et Applications de Metz, Universit{\'e} Paul Verlaine-Metz, B{\^a}t. A, Ile de Saulcy, F-57045 Metz Cedex 01, e-mail: rabaoui@math.univ-metz.fr.}}

\begin{abstract}
In this paper, we prove a generalisation of Bochner-Godement theorem. Our result deals with Olshanski spherical pairs $(G, K)$ defined as inductive limits of increasing sequences of Gelfand pairs $(G(n), K(n))_{n\geq1}$. By using the integral representation theory of G. Choquet on convex cones, we establish a Bochner type representation of any element $\varphi$ of the set ${\mathcal{P}}^{\natural}(G)$ of $K$-biinvariant continuous functions of positive type on $G$. 
\end{abstract}


\maketitle

\section{Introduction}One of the main problems in harmonic analysis is to decompose a unitary representation by means of irreducible ones. The classical Bochner theorem provides an answer for this problem by giving a decomposition of a continuous function of positive type on $\R$ as an integral of indecomposable ones.

 In harmonic analysis on groups of the type $G=\bigcup_{n=1}^\infty
 G(n),$ where $G(n)$ is a sequence of classical groups, with a
 subgroup $K$ of the same type $K=\bigcup_{n=1}^\infty K(n),\
 K(n)\subset G(n),$ several extensions of the Bochner theorem had
 been proved. For example, E. Thoma in 1964 and S. Kerov, G. Olshanski
 and A. Vershik in 2004 studied the case of the infinite symmetric
 group $\mathfrak{S}_\infty= \bigcup_{n=1}^\infty \mathfrak{S}_n$,
 with $G=\mathfrak{S}_\infty\times\mathfrak{S}_\infty$ and $K={\rm diag}(\mathfrak{S}_\infty\times\mathfrak{S}_\infty)\simeq\mathfrak{S}_\infty$ (cf. \cite{Tho2}, \cite{Ol5}). D. Voiculescu in 1976 and G. Olshanski in 2003 treated the pair $G=U(\infty)\times U(\infty)$, $K={\rm diag}(U(\infty)\times U(\infty))\simeq U(\infty)$, where $U(\infty)=\bigcup_{n=1}^\infty U(n)$ is the infinite dimensional unitary group (cf. \cite{Ol4}, \cite{Voi}).

 G. Olshanski proved that the inductive limit of an increasing
 sequence of Gelfand pairs is a spherical pair. Hence, the cited
 examples and many others are part of G. Olshanski's theory for
 spherical pairs which was elaborated in 1990
 (cf. \cite{Ol1}). However, a Bochner type decomposition in this
 setting has not been established yet. In this paper, by using Choquet's theorem, we
 prove such
 generalisation, answering a question asked by J. Faraut in \it
 Infinite Dimensional Harmonic Analysis and Probability \rm (cf. \cite{Far2}).

 This paper consists of 4 sections devoted to the following topics : in section 2 we begin by recalling some definitions and results concerning continuous functions of positive type, then we prove that, for a classical Gelfand pair $(H, M)$, the commutant $\pphi(H)^{'}$ is commutative and use this to give a direct proof of the fact that the set $\P^\bi(H)$ of $M$-biinvariant continuous functions of positive type on $H$ is a lattice. In section 3, we move to the general setting of an increasing sequence of Gelfand pairs $(G(n), K(n))_{n\geq1}$. Our main tool for establishing the generalised Bochner type decomposition is Choquet's theorem. In order to prove the existence of the decomposition, we embed $\P^\bi(G)$, for $G=\bigcup_{n=1}^\infty G(n),$ and $K=\bigcup_{n=1}^\infty K(n),$ into a bigger set $\Q$. For the uniqueness, we prove that the commutant $\pphi(G)^{'}$ remains commutative, and that $\P^\bi(G)$ is a lattice too. At the end of this paper, we present some remarks and open questions.

 We have tried to keep notations and proofs to a minimum in order to
 make the presentation as clear as possible, we refer to \cite{Bo},
 \cite{Go1}, \cite{Go2} and \cite{Go3} for more details on functions
 of positive type and Bochner theorem. The method we follow in our
 proof is a generalisation of E. Thoma's method in the case of a
 countable discrete group (cf. \cite{Tho1}), with some modifications
 inspired from Olshanski's work on the space of infinite dimensional
 hermitian matrices (cf. \cite{Ol2}).
\section{Definitions and results for continuous functions of positive type}
We first recall some definitions and results about functions of positive type. Let $G$ be a Hausdorff topological group having $e$ as unit, and $K$ a closed subgroup of $G$.
\begin{defi}\rm{A function $\pp:G\longrightarrow \C$ is said to be {\it of positive type} if the kernel defined on $G\times G$ by $(g_1,g_2)\longmapsto\pp(g_2^{-1}g_1)$ is of positive type, i.e. for all $g_1, g_2,\ldots, g_n \in G$ and all $c_1, c_2,\dots, c_n \in \C$, $$\sum_{i=1}^n\sum_{j=1}^nc_i\ove{c_j}\pp(g_j^{-1}g_i)\ \geq 0 .$$}
\end{defi}
\begin{pro} Every function $\pp$ of positive type on $G$ is hermitian, i.e. for all $g\in G$, $\ove{\pp(g)}=\pp(g^{-1})$. In addition, $\pp$ is bounded : $|\pp(g)|\leq\pp(e)$.
\end{pro}

A function $\pp$ defined on $G$ is said to be {\it $K$-biinvariant} if it verifies $\pp(k_1gk_2)=\pp(g)$, for all $k_1$, $k_2\in K$ and all $g\in G$. For a unitary representation $(\pii,\H)$, we denote by $\H_K$ the subspace of $K$-invariant vectors in $\H$.
\begin{pro} \label{Pone} Let $(\pii,\H)$ be a unitary representation of $G$ and $\xi$ a vector in $\H_K$. Then, the function $\pp:G\longrightarrow 
\C$ , $g\longmapsto \<\pii(g)\xi,\xi\>_\H$ is $K$-biinvariant of positive type.

\end{pro}

Using the G.N.S. (Gelfand-Naimark-Segal) construction, we can prove that every $K$-biinvariant function of positive type on $G$ can be represented by a unitary representation on $G$.

\begin{pro}\label{Tone}{\bfseries{\rm(}G.N.S. construction{\rm)}} Let
  $\pp$ be a $K$-biinvariant continuous function of positive type on $G$. Then, there exists a triplet $(\pphi,\hphi,\xi^\pp)$ consisting of a unitary representation $\pphi$ on a Hilbert space $(\hphi,\<.,.\>_\pp)$, and a cyclic vector $\xi^\pp\in\hphik$ such that, for all $g\in G$, $$\pp(g)=\<\pphi(g)\xi^\pp,\xi^\pp\>_\pp.$$ Moreover, this triplet is unique in the following sense : if 
$(\pii,\H,\xi)$ is another triplet, then there exists an interwining isomorphism $T : \hphi\rightarrow\H$ between $\pphi$ and $\pii$ such that $T\xi^\pp=\xi$.
\end{pro}

Let $\P(G)$ be the set of continuous functions of positive type on $G$. $\P(G)$ is a convex cone which is invariant under product and complex conjugation.

For a convex set $E$, we denote by Ext($E$) its subset of extremal
points. We also denote by $\P_{\leq 
1}(G)$ (respectively $\P_1(G)$) the set of elements $\pp$ of $\P(G)$ verifying $\pp(e)\leq 1$ (respectively $\pp(e)=1$).

\begin{lem}\label{Lone}\rm{Ext(}$\P_{\leq 1}(G)$\rm{) = Ext(}$\P_1(G)$\rm{)} $\cup$ \rm{\{}$0$\rm{\}}.
\end{lem}

Next, we will prove some algebraic characterizations which will be used to establish the uniqueness of the decomposition given by the generalized Bochner theorem.

Let $\Gamma$ be a convex cone in a topological vector space $E$. This cone is equipped with its proper order : $\gamma_1\ll\gamma_2$ if $\gamma_2-\gamma_1\in\Gamma$. The cone $\Gamma$ is said to be a {\it lattice} if each couple of elements $\gamma_1$, $\gamma_2$ in $\Gamma$ have (for the order defined by the cone) {\it a least upper bound} in $\Gamma$, denoted by $\gamma_1\vee\gamma_2$, and {\it a greatest lower bound} in $\Gamma$, denoted by $\gamma_1\wedge\gamma_2$.\\ For $\gamma_0\in\Gamma$, we denote by $\Gamma^{\gamma_0}$ {\it the face of $\Gamma$} defined as: $$\Gamma^{\gamma_0}=\{\gamma\in\Gamma\  |\  \exists\ \lambda\geq 0\  ;\ 
\gamma\ll\lambda\gamma_0\}.$$ The order of $\Gamma^{\gamma_0}$ coincides with the one induced by $\Gamma$. The cone $\Gamma$ is a lattice if and only if, for every $\gamma_0$, the face 
$\Gamma^{\gamma_0}$ is a lattice. 

Let now $\Gamma=\P^\bi(G)$ be the subcone of $\P(G)$ which consists of $K$-biinvariant elements. On this convex cone, and similarly on $\Pin$, the proper order $\ll$ is given by: $$\pp\ll \psi \quad \textrm{if and only if} \quad \psi-\pp\in \P^\bi(G) \quad (\pp,\psi\in\P^\bi(G)).$$ Recall that every function $\pp\in\P^\bi(G)$ is associated to a triplet $(\pphi,\hphi,\xi^\pp)$. Let $\mathcal{A}=\pphi(G)^{'}$ be {\it the commutant} of $\pphi(G)$. It is a selfadjoint subalgebra of $\L(\hphi)$. We will prove that each face $\Gamma^{\pp}$ of $\P^\bi(G)$ is lineary isomorphic to the cone $\A^{+}=\{T\in\A\ |\ \forall\ v\in\hphi, \<Tv,v\>_\pp\geq 0\}$ of positive operators of $\A$ on which we define an order, denoted $\prec$ : $$P\prec Q \quad \textrm{if and only if} \quad  \<Pv,v\>_\pp \ \leq \ \<Qv,v\>_\pp \quad (v\in\hphi,\ P,Q\in\A^+).$$ 

\begin{teo}\label{Theo1} Let $K$ be a closed subgroup of a Hausdorff topological group $G$. For all $\pp\in\P^\bi(G)$ the face $\Gamma^\pp$ is 
lineary isomorphic to the cone $\A^{+}$ of positive operator of the algebra $\A=\pphi(G)^{'}$. This bijective correspondence identifies an element $\psi\in\Gamma^\pp$ with an element 
$T\in\A^{+}$ such that   \bq\label{eqone}\psi(g)=\<T\pphi(g)\xi^\pp,\xi^\pp\>_\pp, \ g\in G.\eq 
\begin{paragraph}{Proof.}\rm{Let $T\in\A^{+}$. The operator
    $T^{\frac{1}{2}}$ exists and belongs to $\A^{+}$ (\cite{Do}, page
    430, 11.17). So, for all $g\in G$,
    \ba\psi(g)\ =\ \<T\pphi(g)\xi^\pp,\xi^\pp\>_\pp&=&\<T^{\frac{1}{2}}\pphi(g)\xi^\pp,(T^{\frac{1}{2}})^*\xi^\pp\>_\pp \\&=&\<\pphi(g)T^{\frac{1}{2}}\xi^\pp,T^{\frac{1}{2}}\xi^\pp\>_\pp .\ea  The function $\psi$ is of positive type (Proposition 2). It is also continuous since the map $\xi\longmapsto\pphi(g)\xi$ is continuous for every $g\in G$. It is also $K$-biinvariant. Hence, $\psi\in\P^\bi(G)$.

 If we put $\lambda_0=||T||$, where $||.||$ is the usual operator norm defined on $\L(\hphi)$, then $\lambda_0\pp-\psi\in\P^\bi(G)$. In fact \ba 
(\lambda_0\pp-\psi)(g)&=&||T||\<\pphi(g)\xi^\pp,\xi^\pp\>_\pp-\<\pphi(g)T\xi^\pp,\xi^\pp\>_\pp\\
&=&\<\pphi(g)C\xi^\pp,\xi^\pp\>_\pp,\ea where $C=||T||I-T$. As, for all $v\in\hphi$, $0\leq \<Tv,v\>_\pp \leq ||T||\<v,v\>_\pp$, the operator $C\in\A^{+}$. Hence $C=D^2$ with $D\in\A^{+}$, and so $$(\lambda_0\pp-\psi)(g)=\<\pphi(g)D^2\xi^\pp,\xi^\pp\>_\pp= 
\<\pphi(g)D\xi^\pp,D\xi^\pp\>_\pp.$$ This proves, by Proposition 2, that $\lambda_0\pp-\psi$ is of positive type. It is also continuous and $K$-biinvariant. Hence, $\lambda_0\pp-\psi\in\P^\bi(G)$.

One can also remark that $\psi$ uniquely determine $T$. In fact, for all $g,h\in G$, $$\psi(h^{-1}g)=\<\pphi(h^{-1}g)T\xi^\pp,\xi^\pp\>_\pp= \<T\pphi(g)\xi^\pp,\pphi(h)\xi^\pp\>_\pp.$$ If 
$\widetilde{T}$ is another operator in $\A^{+}$ verifying (\ref{eqone}), then for all $g,h\in G$, $$\<\pphi(g)(T-\widetilde{T})\xi^\pp,\pphi(h)\xi^\pp\>_\pp=0.$$ Since $V_\pp=Vect\{\pphi(g)\xi^\pp \ , \ g\in G\}$ is dense in $\hphi$, $$T=\widetilde{T}.$$

It remains to prove that, for every $\psi\in\Gamma^\pp$, there exists $T\in\A^{+}$ verifying (\ref{eqone}). Let us denote by $$\mathfrak{M^0}(G):=\{\m=\sum_{i=1}^ma_i\delta_{x_i}\ | \ (a_i)_{i}\subset\C \ ,\ (x_i)_{i}\subset G\},$$ the space of measures with finite support. For a function of positive type $\pp$ and $\m, \n\in\mathfrak{M^0}(G)$, put $$(\pp,\n^**\m)=\sum_{i=1}^m\sum_{j=1}^n\ove{b_j}a_i\pp(x_j^{-1}x_i) \geq \ 0.$$  We can also define the function $$\m*\pp(x)=\into 
\pp(y^{-1}x)d\m(y)=\sum_{i=1}^ma_i\pp(x_i^{-1}x),$$ it is continuous
and right $K$-invariant. With the previous notation and
definitions, the vector space $V_\pp$ can also be given by :
$$V_\pp:=\{\pp^\m=\m*\check{\pp}=\sum_{i=1}^ma_i\pphi(g_i)\xi^\pp, \
\m\in \mathfrak{M^0}(G)\},$$ where $\check{\pp}(g)=\pp(g^{-1})$, for
all $g\in G$. For $\pp^\m, \pp^\n\in V_\pp$, put
$$\<\pp^\m,\pp^\n\>_\pp = (\pp,\n^**\m).$$ The map $(\pp^\m,\pp^\n)\longmapsto \<\pp^\m,\pp^\n\>_\pp$ is a hermitian positive form on $V_\pp$, which is in addition definite as it verifies, for all $g\in G$, $$|\pp^\m(g)|^2=|\m*\pp(g)|^2\leq \pp(e)\<\pp^\m,\pp^\m\>_\pp.$$ Now, let $\psi\in\Gamma^\pp$, there exists $\lambda_0\geq 0$ such that $$\lambda_0\pp-\psi\in\P^\bi(G).$$ So, for all $\m\in \mathfrak{M^0}(G)$,  $$(\lambda_0\pp 
-\psi,\m^**\m) \geq 0 \ \ \rm{or \ equivalently} \ \ 
(\psi,\m^**\m)\leq(\pp,\m^**\m).$$ Hence $$\<\psi^\m,\psi^\m\>_\psi \ \leq 
\lambda_0\<\pp^\m,\pp^\m\>_\pp.$$ Consequently, we can define on $V_\pp\times V_\pp$ a hermitian form $\omega$ given, for every $\m, \n\in\mathfrak{M^0}(G)$, by $$\omega(\pp^\m,\pp^\n)=(\psi,\n^**\m)=\<\psi^\m,\psi^\n\>_\psi.$$ In fact $$|\omega(\pp^\m,\pp^\n)|^2=|\<\psi^\m,\psi^\n\>_\psi|^2\leq\lambda_0^2\<\pp^\m,\pp^\m\>_\pp\<\pp^\n,\pp^\n\>_\pp.$$ In addition 
$$\omega(\pp^\m,\pp^\n)=(\psi,\n^**\m)= 
\ove{(\psi,\m^**\n)}=\ove{\omega(\pp^\n,\pp^\m)}.$$ So, $\omega$ is a well-defined hermitian form which is continuous on $V_\pp\times V_\pp$. It is also positive as, for all $\m\in \mathfrak{M^0}(G)$, $$\omega(\pp^\m,\pp^\m)=(\psi,\m^**\m)\geq 
0.$$ 

As $V_\pp$ is dense in $\hphi$, $\omega$ may be extended to a positive hermitian continuous form on $\hphi\times\hphi$. So, by Riesz's theorem, there exists an unique positive hermitian operator $T$ in $\L(\hphi)$ such that, for every $v_1, v_2\in\hphi$, $$\<Tv_1,v_2\>_\pp=\omega(v_1,v_2).$$ In particular, for $\pp^\m,\pp^\n\in V_\pp$, $$\<T\pp^\m,\pp^\n\>_\pp=\omega(\pp^\m,\pp^\n)=(\psi,\n^**\m).$$ Consequently, for $\m_0=\delta_g\ ,\  g\in G$ and $\n_0=\delta_e$, $$\<T\pp^{\m_0},\pp^{\n_0}\>_\pp=\<T\pp^{\delta_g},\pp^{\delta_e}\>_\pp=(\psi,\delta_e^**\delta_g)=\psi(g).$$ But, $\pp^{\delta_g}=\pphi(g)\xi^\pp$ and $\pp^{\delta_e}=\xi^\pp$. Hence $\psi(g)=\<T\pphi(g)\xi^\pp,\xi^\pp\>_\pp.$ The operator $T$ is also selfadjoint and positive. In fact, as $\psi$ is of positive type, for every $g, h\in G$, $\psi(g^{-1}h)=\ove{\psi(h^{-1}g)}$. Hence $$\<T\pphi(h)\xi^\pp,\pphi(g)\xi^\pp\>_\pp=\ove{\<T\pphi(g),\pphi(h)\xi^\pp\>_\pp },$$ and so 
$$\<\pphi(h)\xi^\pp,T^*\pphi(g)\xi^\pp\>_\pp=\<\pphi(h)\xi^\pp,T\pphi(g)\>_\pp .$$ Since $V_\pp$ is dense in $\hphi$, $$T=T^*.$$ The positivity of $T$ follows from $\omega$'s one. The operator $T$ also commutes with $\pphi(g)$, for all $g\in G$. $\hspace{1cm} \Box$}                                                        
\end{paragraph}
\end{teo}
Next, we give a necessary and sufficient condition for the cone $\P^\bi(G)$ to be a lattice.
\begin{lem}The cone $\A^{+}$ is a lattice if and only if the algebra $\A$ is commutative.
\begin{paragraph}{Proof.}\rm{The proof is similar to the one given in (\cite{Far3}, Theorem III.2.4, page 129). $\hspace{1cm}\Box$    } 
\end{paragraph}
\end{lem}
By Theorem \ref{Theo1} and this last lemma, we prove the following theorem,
\begin{teo}\label{Theo2}Let $K$ be a closed subgroup of a Hausdorff topological group $G$. The cone $\P^\bi(G)$ is a lattice if and only if, for every function $\pp$ of this cone, the algebra $\A=\pphi(G)^{'}$ is commutative.
\begin{paragraph}{Proof.} \rm{From Theorem \ref{Theo1}, we deduce that, for every function $\pp\in\P^\bi(G)$, the face $\Gamma^\pp$ is lineary isomorphic to the cone $\A^{+}$, which is a lattice if and only if $\A$ is commutative. So, for every function $\pp\in\P^\bi(G)$, $\Gamma^\pp$ is a lattice if and only if $\A$ is commutative. $\hspace{1cm}\Box$ }
\end{paragraph}
\end{teo}
\begin{defi}\rm{A pair $(G,K)$, where $G$ is a locally compact group and $K$ a compact subgroup of $G$, is said to be a {\it Gelfand pair} if the convolution algebra of $K$-biinvariant integrable functions is commutative. }
\end{defi}
 We will prove by using some elements of
von Neumann algebra theory that, in the case of a Gelfand pair $(G,K)$, the algebra $\pphi(G)^{'}$ is commutative, for all $\pp\in\P^\bi(G)$.
\begin{pro}Let $(G,K)$ be a Gelfand pair and $P$ the orthogonal projection onto $\hphi_K$ defined by $$P=\int_K\pphi(k)\ \alpha(dk),$$ where $\alpha$ is the normalized Haar measure of the subgroup $K$. Then $P$ is an element of $\pphi(G)^{''}$, and the algebra $P\pphi(G)^{''}P$ is commutative.
\begin{paragraph}{Proof.}\rm{ Let us prove that $P\in\pphi(G)^{''}$. In fact, for every $T\in\pphi(G)^{'}$ and every $v,w\in\hphi$,  $$ \<PTv,w\>=\<\pphi(\alpha)Tv,w\>=\<\pphi(\alpha)v,T^*w\>=\<TPv,w\>.$$ So, for every $v$ in $\hphi$, $PTv=TPv$. Hence $P\in\pphi(G)^{''}$. As $(G,K)$ is a Gelfand pair, for every $\mu,\ \nu \in \mathfrak{M^0}(G)$, the $K$-biinvariant measures $\alpha*\mu*\alpha$ and $\alpha*\nu*\alpha$ commute. Thus, for every $\mu,\ \nu \in \mathfrak{M^0}(G)$, 
$$P\pphi(\mu)P\pphi(\nu)P=P\pphi(\nu)P\pphi(\mu)P.$$ As $\pphi(\mathfrak{M^0}(G))$ is a selfadjoint subalgebra containing the identity of $\L(\H^\pp)$, it is dense in $\pphi(G)^{''}$ in the strong topology of operators (\cite{D1}, Theorem 2 and Corollary 1, page 45). Hence, for every $A,\ B\in\pphi(G)^{''}$, $$PAPBP=PBPAP.$$ Put $S=PAP$ and $T=PBP$. $S$ and $T$ are two arbitrary elements of the algebra $P\pphi(G)^{''}P$ and they verify 
$$ST=PAPPBP=PAPBP=PBPAP=TS.$$ It follows that the algebra $P\pphi(G)^{''}P$ is commutative.$\hspace{1cm}\Box$}
\end{paragraph}
\end{pro}

For an operator $A$ of the von Neumann algebra $\pphi(G)^{'}$, let us denote by $A_P$ the restriction of the operator $PA$ to $\hphi_K$. Put $[\pphi(G)^{'}]_P=\{A_P, \ A\in\pphi(G)^{'}\}$. By (\cite{D1}, Proposition 1, page 18), the algebras $[\pphi(G)^{'}]_P$ and $[\pphi(G)^{''}]_P$ are von Neumann algebras and they verify $$([\pphi(G)^{''}]_P)^{'}=[\pphi(G)^{'}]_P.$$ 

Since $\xi^\pp$ is a cyclic vector for the algebra $\pphi(\mathfrak{M^0}(G))$, by (\cite{D2}, Appendice A, A14), it is a separating vector for the von Neumann algebra $\pphi(\mathfrak{M^0}(G))^{'}=\pphi(G)^{'}$. Thus it is also separating for the von Neumann algebra $[\pphi(G)^{'}]_P$. Hence it is cyclic for the von Neumann algebra $[\pphi(G)^{''}]_P$. 

By using the fact that every von Neumann algebra $\M$ which is commutative and possesses a cyclic vector verifies $\M^{'}=\M$ (\cite{D1}, Corollaire 2, page 89), and by noticing that the algebra $[\pphi(G)^{''}]_P$ is nothing but $P\pphi(G)^{''}P$, we obtain $$([\pphi(G)^{''}]_P)^{'}=[\pphi(G)^{''}]_P.$$ Hence $$[\pphi(G)^{'}]_P=[\pphi(G)^{''}]_P.$$  Now, to get the commutativity of $\pphi(G)^{'}$, it is sufficient to prove the following proposition,
\begin{pro}Let $(G,K)$ be a Gelfand pair. The commutant $\pphi(G)^{'}$, seen as a von Neumann algebra, is isomorphic to the algebra $[\pphi(G)^{'}]_P$.
\begin{paragraph}{Proof.}\rm{Let $\Psi : \pphi(G)^{'}\rightarrow[\pphi(G)^{'}]_P$, $A\longmapsto A_P$. $\Psi$ is well-defined, it is also a homomorphism of algebras, since for every $S,\ T\in\pphi(G)^{'}$,  $$\Psi(ST)=[ST]_P=PSTP=PSPPTP=S_PT_P=\Psi(S)\Psi(T).$$ And
    $$\Psi(T^*)=PT^*P=P^*T^*P^*=(PTP)^*=(T_P)^*=\Psi(T)^*.$$ It is evident that $\Psi$ is onto by construction. Let us prove that it is one to one. 

Let $S\in\pphi(G)^{'}$ such that $\Psi(S)=0$. Then, $$\Psi(S)=0 \Rightarrow PS\xi^\pp=0\Rightarrow 
SP\xi^\pp=0\Rightarrow S\xi^\pp=0.$$ Hence, for every $g\in G$, $S\pphi(g)\xi^\pp=\pphi(g)S\xi^\pp=0.$ And since $\xi^\pp$ is cyclic, we get immediately $S=0$. Therefore, $\Psi$ is one to one.$\hspace{1cm}\Box$  }
\end{paragraph}
\end{pro}

\begin{teo}\label{Theo3}Let $(G,K)$ be a Gelfand pair and $\pp$ a $K$-biinvariant continuous function of positive type on $G$. Then, the algebra $\pphi(G)^{'}$ is commutative.
\begin{paragraph}{Proof.}\rm{By the previous proposition, $\pphi(G)^{'}$ is isomporphic to $[\pphi(G)^{'}]_P$. Also we know that $[\pphi(G)^{'}]_P=[\pphi(G)^{''}]_P=P\pphi(G)^{''}P.$ The result follows since the algebra $P\pphi(G)^{''}P$ is commutative.$\hspace{1cm}\Box$}
\end{paragraph}
\end{teo}

\begin{cor}\label{Cor1} Let $(G,K)$ be a Gelfand pair. Then, the cone $\P^\bi(G)$ is a lattice.
\begin{paragraph}{Proof.}\rm{By Theorem \ref{Theo2}, $\P^\bi(G)$ is a lattice if and only if, for every element $\pp$ in this cone, the algebra $\pphi(G)^{'}$ is commutative, which is satisfied in this case as shown by the previous theorem. Hence $\P^\bi(G)$ is a lattice.$\hspace{1cm}\Box$}
\end{paragraph}
\end{cor}
We know that every function of positive type is bounded. Since $G$ is
a locally compact topological group, $\P(G)$ can be seen as a subset
of $L^\infty(G)$ for a left invariant Haar measure on $G$. We add,
from now on, the condition that $G$ is separable and we consider on
$\P(G)$ the topology induced by the weak-$*$ topology $\sigma
(L^\infty(G),L^1(G))$, denoted by $\tau^*(L^\infty(G))$. By the
Banach-Alaoglu theorem (cf. \cite{Rud1}), the unit ball of
$L^\infty(G)$ is compact in this topology. In addition, $\Pin$
considered as a subset of $L^\infty(G)$, is closed in this
topology(cf. \cite{Rud1}, \cite{Far1}). Therefore, $\Pin$ is
compact. Furthermore, the unit ball of $L^\infty(G)$, for $G$
separable, is metrisable in the weak-$*$ topology
$\tau^*(L^\infty(G))$ (cf. \cite{D2}, \cite{Rud1}). Hence $\Pin$ is metrisable.

\newpage

Thus $\Pin$ is convex, compact and metrisable in the topological space $L^\infty(G)$ which is locally convex in the weak-$*$ topology $\tau^*(L^\infty(G))$. Furthermore, by Corollary 1, the cone generated by $\Pin$, namely $\P^\natural(G)$, is a lattice. Therefore, we get by applying Choquet's theorem that every element $\pp\in\P^\natural(G)$ has an integral representation : $$\pp(g)=\int_{\rm{Ext(}\Pun\rm{)}}\omega(g)\m(d\omega).$$ This last formula constitutes Bochner-Godement's theorem. It is evident now that Choquet's theorem is fundamental for the proof. Because of its importance, we finish this section by giving its statement.

\begin{teo}{{\rm(}{\bfseries Choquet's theorem}, see \cite{Phel}
    sections 3 and 10{\rm)}} Let $\mathcal{U}$ be a convex subset of a
  locally convex topological vector space $E$. If $\mathcal{U}$ is
  compact and metrisable, then 

\item[ {\rm (i)}] {\rm Ext(}$\mathcal{U}${\rm)} is a Borel subset of $\mathcal{U}$.

\item[ {\rm (ii)}] For every $a\in \mathcal{U}$, there exists a probability measure $\mu$ on {\rm Ext(}$\mathcal{U}${\rm)}, such that for all continuous linear form $L$ on $E$, $$L(a) = \int_{b\in {\rm Ext(}\mathcal{U}{\rm)}}  
L(b) \mu(db).$$

\item[ {\rm (iii)}] $\m$ is unique if and only if the cone generated by $\mathcal{U}$ is a lattice.
\end{teo}

\section{A Bochner type theorem for Olshanski spherical pairs}

\begin{defi}\rm{Let $H$ be a Hausdorff topological group and $M$ a
    closed subgroup of $H$. The pair $(H,M)$ is said to be {\it spherical } if, for every irreducible unitary representation $\pii$ of $H$ on a Hilbert space $\H$, $${\rm dim}\ \H_M\leq1.$$ If $H$ is locally compact, and $M$ compact, then the pair $(H,M)$ is spherical if and only if it is a Gelfand pair.}
\end{defi}

Let $\big(G(n),K(n)\big)_{n\geq 1}$ be a sequence of Gelfand pairs such that $G(n)$ is a locally compact topological group which is also a closed subgroup of $G(n+1)$. $K(n)$ is a closed subgroup of $K(n+1)$ verifying $K(n)=K(n+1)\cap G(n)$. The family of Gelfand pairs $\big(G(n),K(n)\big)_{n\geq 1}$, equiped with the system of canonical
continuous embeddings from $G(i)$ to $G(j)$ with $ i\leq j$ ,
constitute an inductive countable system of topological groups
(cf. \cite{Bou1}). Hence we may define the following inductive limit
groups : $G=\bigcup_{n=1}^\infty G(n)$ and $K=\bigcup_{n=1}^\infty
K(n)$. The topology defined on $G$ is the inductive limit topology. It
is the finest topology such that all the canonical embeddings from
$G(n)$ into $G$ are continuous. Olshanski proved that $(G,K)$ is a
spherical pair (cf. \cite{Far2}, \cite{Ol1}). Hence we can introduce the following definition: 

\begin{defi} Let $\big(G(n),K(n)\big)_{n\geq 1}$ be an increasing sequence of Gelfand pairs as above. The inductive limit pair $(G,K)$ is called an {\it Olshanski spherical pair.}
\end{defi}
The group $G$ equipped with the inductive limit topology is Hausdorff. But, such topology does not make $G$ locally compact. Therefore we can not directly apply Choquet's theorem to $\P^\bi(G)$ as in the classical case. In order to solve this problem, we embed $\P^\bi(G)$ in the cone of subprojective systems : $$\mathcal{Q} := 
\left\{\pp=\{\pp^{(i)}\}_i\in \prod_{i=1}^\infty\P^\bi(G(i)) \ | \ \res_i^{i+1}\big(\pp^{(i+1)}\big)\ll \pp^{(i)} \ i=1,2,...\right\},$$ where $\res_n^{n+1}$ is the restriction to $G(n)$ of a function defined on $G(n+1)$. Choquet's theory of integral representation applied to $\mathcal{Q}$ will give us a Bochner type theorem for the spherical pairs of Olshanski. Let $\res_n$ be the restriction to $G(n)$ of a function defined on $G$, and put $\P_m^n= 
\prod_{k=m}^n \P^\bi(G(k))$, where $1\leq m \leq n \leq \infty$.

\begin{nota}\rm{If $G_1\subset G_2$ are two locally compact groups the set of pairs $\{(\pp, \psi)\in\P(G_1)\times\P(G_2)\ | \ \res(\psi) = \pp\}$, where $\res$ is the restriction to $G_1$ of a function on $G_2$, is not closed in general, and in some cases it can be shown that it is dense in $\{(\pp, \psi)\in\P(G_1)\times\P(G_2)\ | \ \res(\psi) \ll \pp\}$.}
\end{nota}

Next we will prove that $\Q$ is closed in $\P_1^\infty$ in the product
topology $\tau^*=\prod_{n=1}^\infty\tau^*(L^\infty(G(n)))$. To
establish this, it is sufficient to prove that the set
$$\mathcal{R}_k=\left\{(\pp^{(k)},\pp^{(k+1)})\in\P_k^{k+1} \ |\
  \res_k^{k+1}(\pp^{(k+1)})\ll\pp^{(k)}\right\}$$ is closed in the
topology $\tau^*(L^\infty(G(k)))\times\tau^*(L^\infty(G(k+1)))$.

 Let $H$ be a locally compact group, $\alpha$ its left invariant Haar measure, and $M$ a compact subgroup of $H$ such that $(H,M)$ is a Gelfand pair.
\begin{lem}\label{Lem3} For every function $\pp\in\P^\bi(H)$ and $f\in L^1(H)^\bi$ such that $||f||_1\leq 1$, one has $$f^**\pp*f\ll\pp.$$
\begin{paragraph}{Proof.}\rm{Let $(\pphi, \hphi)$ be the unitary representation associated to $\pp$ : $$\pp(h)=\<\pphi(h)\xi^\pp, \xi^\pp\>_\pp \ \ \ (h\in H).$$ Since $(H,M)$ is a Gelfand pair, the operator $\pphi(f)$ commutes, for every $h\in H$, with $\pphi(h)$, and $$f^**\pp*f(h)=\<\pphi(h)\pphi(f)\xi^\pp, \pphi(f)\xi^\pp\>_\pp.$$ Therefore \ba\sum_{i,j=1}^Nf^**\pp*f(h_j^{-1}h_i)c_i\ove{c_j}&=&||\sum_{i=1}^Nc_i\pphi(h_i)\pphi(f)\xi^\pp||_\pp^2\\
                                                   &=&||\pphi(f)\sum_{i=1}^Nc_i\pphi(h_i)\xi^\pp||_\pp^2\\
                                                   &\leq&||\pphi(f)||^2||\sum_{i=1}^Nc_i\pphi(h_i)\xi^\pp||_\pp^2\\
                                                   &\leq&||\sum_{i=1}^Nc_i\pp(h_i)\xi^\pp||_\pp^2\\
                                                   &=&\sum_{i,j=1}^N\pp(h_j^{-1}h_i)c_i\ove{c_j}. \hspace{1cm}\Box \ea}
\end{paragraph}
\end{lem}
Under the same assumptions as Lemma \ref{Lem3}, we prove the following lemma,
\begin{lem}\label{Lem4}The linear form $L$ defined, for every bounded measure $\m$ on $H$, by $$L(\pp)=\int_{H\times H}\pp(y^{-1}x)\m(dx)\ove{\m(dy)}$$ is lower-semicontinuous on $\P^\bi(H)$ in the weak-$*$ topology $\tau^*(L^\infty(H))$.
\begin{paragraph}{Proof.}\rm{Firstly, let us remark that $L$ is positive on $\P^\bi(H)$ and that if $\m=\delta$, then $L(\pp)=\pp(e)$. We will prove that, for every constant $C\geq 0$, the set $$\{\pp\in\P^\bi(H)\ | \ L(\pp)\leq C\}$$ is closed. Let $(\pp_n)$ be a sequence of $\P^\bi(H)$ that converges to $\pp$, i.e. for every $f\in L^1(H)$, $$\lim_{n\rightarrow\infty}\int_H\pp_n(h)f(h)\alpha(dh)=\int_H\pp(h)f(h)\alpha(dh).$$ Suppose that, for every $n$, $L(\pp_n)\leq C$. We know that, for every bounded measure $\m$ and $f\in L^1(H)^\bi$, $f*\m\in L^1(H)$. Suppose $||f||_1\leq 1$. By hypothesis, for every $n$, $$\m^**\pp_n*\m(e)\leq C.$$ Therefore, by Lemma \ref{Lem3}, $$\m^**f^**\pp_n*f*\m(e)\leq C,$$ and since $$\lim_{n\rightarrow\infty}\m^**f^**\pp_n*f*\m(e)=\m^**f^**\pp*f*\m(e),$$ it follows that $$\m^**f^**\pp*f*\m(e)\leq C.$$ By considering an approximation of the identity $(f_k)$ : $f_k\in L^1(H)^\bi$, $f_k\geq 0$,$$\int_Hf_k(h)\alpha(dh)=1,$$ and observing that for every continuous bounded function $\psi$ : $$\lim_{k\rightarrow\infty}\int_H\psi(h)f_k(h)\alpha(dh)=\psi(e),$$ we deduce that $$\m^**\pp*\m(e)\leq C.\hspace{1cm}\Box$$ }
\end{paragraph}
\end{lem}

\begin{pro}Let $U$ be a closed unimodular subgroup of $H$, $\alpha_U$ its left invariant Haar measure and $\res$ the application that for a function on $H$ associates its restriction to $U$. The set $$\{(\phi,\psi)\in\P^\bi(H)\times\P^\bi(U)\ |\ \res(\phi)\ll\psi\}$$ is closed.
\begin{paragraph}{Proof.}\rm{Let $(\phi_n,\psi_n)$ be a sequence in $\P^\bi(H)\times\P^\bi(U)$ that converges to $(\phi,\psi)$, and suppose that, for every $n$ and every function $f\in L^1(U)$, $$ \int_{U\times U}\phi_n(y^{-1}x)f(x)\ove{f(y)}\alpha_U(dx)\alpha_U(dy) \leq \int_{U\times U}\psi_n(y^{-1}x)f(x)\ove{f(y)}\alpha_U(dx)\alpha_U(dy).$$ Let $$C > \int_{U\times U}\psi(y^{-1}x)f(x)\ove{f(y)}\alpha_U(dx)\alpha_U(dy).$$ There exists $n_0$ such that, if $n\geq n_0$ $$\int_{U\times U}\psi_n(y^{-1}x)f(x)\ove{f(y)}\alpha_U(dx)\alpha_U(dy)\leq C,$$ and thus $$\int_{U\times U}\phi_n(y^{-1}x)f(x)\ove{f(y)}\alpha_U(dx)\alpha_U(dy) \leq C.$$ Lemma \ref{Lem4} applied to the measure $\m(dx)=f(x)\alpha_U(dx)$ gives $$\int_{U\times U}\phi(y^{-1}x)f(x)\ove{f(y)}\alpha_U(dx)\alpha_U(dy)\leq C.$$ This being true for every constant $C$ verifying  $$C > \int_{U\times U}\psi(y^{-1}x)f(x)\ove{f(y)}\alpha_U(dx)\alpha_U(dy),$$ we can deduce that $$\int_{U\times U}\phi(y^{-1}x)f(x)\ove{f(y)}\alpha_U(dx)\alpha_U(dy) \leq \int_{U\times U}\psi(y^{-1}x)f(x)\ove{f(y)}\alpha_U(dx)\alpha_U(dy).$$ Therefore $\res(\phi)\ll\psi$. It follows that the set $$\{(\phi,\psi)\in\P^\bi(H)\times\P^\bi(U)\ |\ \res(\phi)\ll\psi\}$$ is closed. $\hspace{1cm}$ $\Box$\\}
\end{paragraph}
\end{pro}

Since, for all $n$, the pair $(G(n),K(n))$ is supposed to be a Gelfand pair, the groups $G(n)$ are all unimodular (see \cite{Far1}, Proposition I.1). Hence we can apply the previous proposition in the case where $H=G(k+1)$ and $U=G(k)$. Then, one gets that $\mathcal{R}_k$ is closed, for every $k$, and hence $\Q$ is closed in $\P_1^\infty$. As a consequence, the set $$\mathcal{Q}_{\leq 1} := 
\left\{\pp=\{\pp^{(i)}\}_i\in \prod_{i=1}^\infty\P_{\leq 1}^\bi(G(i)) \ | \ \res_i^{i+1}\big(\pp^{(i+1)}\big)\ll \pp^{(i)} \ i=1,2,...\right\},$$ is compact. In order to get the metrisability of $\Q_{\leq1}$, it is sufficient to suppose that all the $G(n)$ are separable.\\

 It remains to prove that the cone $\Q$ is a lattice in order to apply Choquet's theorem.\\

Let $(\pphi,\hphi,\xi^\pp)$ be the triplet associated to a function
$\pp\in\P^\bi(G)$. We are going to prove that the algebra
$\pphi(G)^{'}$ is commutative. Since $G(n)$ is a subgroup of $G$, the representation $\pphi$ of $G$ remains
a continuous unitary representation of $G(n)$ on $\hphi$. Put
$\displaystyle{\hphi_n=\ove{Vect\{\pphi(g)\xi^\pp \ ,\ g\in G(n)\}}}$. It is a $G(n)$-invariant closed subspace of $\hphi$. Hence we may restrict, for every $g\in G(n)$, the operator $\pphi(g)$ to $\hphi_n$. We obtain a continuous unitary representation of $G(n)$ on $\hphi_n$ that will be denoted by $\pphi_n$.\\ 

 Let $P_n$ be the orthogonal projection onto $\hphi_n$,
\begin{lem}
\item[ {\rm (i)}] $\displaystyle{\bigcup_{n=1}^\infty\hphi_n}$ is dense in $\hphi$.
\item[ {\rm (ii)}] $P_n$ converges strongly to the identity $I$ of $\hphi$.
\end{lem}

\begin{pro} Let $(G,K)$ be an Olshanski spherical pair. For every $\pp\in\P^\bi(G)$, the commutant $\A=\pphi(G)^{'}$ of the representation $\pphi$ which is associated to $\pp$ by the G.N.S. construction, is a commutative algebra.
\begin{paragraph}{Proof.}\rm{Let $B$ be an arbitrary operator of $\A$. Then, for every $g$ in $G$, $B$ commutes with $\pphi(g)$. This is also true on $G(n)$, for every $n\in\mathbb{N}^*$. On the other hand, for every $n\in\mathbb{N}^*$, $P_nBP_n$ which is an operator of $\L(\hphi_n)$ commutes with the representation $\pphi_n$ of $G(n)$ on $\hphi_n$. 
Since $\hphi_n$ is $G(n)$-invariant, for every $g\in G(n)$, $P_n$ commutes with $\pphi(g)$. Therefore, for every $g\in G(n)$, $$P_nBP_n\pphi_n(g)=P_nB\pphi_n(g)P_n=P_n\pphi_n(g)BP_n=\pphi_n(g)P_nBP_n.$$
    By Theorem \ref{Theo3}, the algebra $\pphi_n(G(n))^{'}$ is commutative. So, for two operators $B_1$ and $B_2$ of $\pphi(G)^{'}$, and for every $n\in\mathbb{N}^*$, $$P_nB_1P_nP_nB_2P_n=P_nB_2P_nP_nB_1P_n,$$  $$
    P_nB_1P_nB_2P_n=P_nB_2P_nB_1P_n.$$ Since $K_n\subset
    K_{n+1}$, then $\H_{K_{n+1}}\subset\H_{K_n}$, and therefore
    $$P_{n+1}=P_nP_{n+1}=P_{n+1}P_n.$$ Also, for every $n,m\geq 1$, $$P_{n+m}=P_{n+m}P_n=P_nP_{n+m}.$$ Hence, for every $m, m', n\geq 1$, $$P_{n+m}B_1P_nB_2P_{n+m'}=P_{n+m}B_2P_nB_1P_{n+m'}.$$
By using the fact that $P_n$ converges strongly to $I$ and by pushing $m$, $m'$ to $\infty$, one obtains
$$B_1P_nB_2=B_2P_nB_1.$$ Finally, by pushing $n$ to $\infty$, one gets 
$$B_1B_2=B_2B_1.\hspace{1cm}\Box $$       }
\end{paragraph} 
\end{pro}
\begin{teo}\label{Theo5}For an Olshanski spherical pair $(G, K)$, the cone $\P^\bi(G)$ is a lattice.
\begin{paragraph}{Proof.}\rm{By the previous proposition, the algebra
    $\A=\pphi(G)^{'}$ is commutative. Hence, by Theorem \ref{Theo2}, the cone $\P^\bi(G)$ is a
    lattice.$\hspace{1cm}\Box$}
\end{paragraph}
\end{teo}
Let us prove that $\Q$ is a lattice. We start by giving a decomposition of the elements of $\Q$.

\begin{lem}\label{Lem6} Let $H$ be a locally compact topological group having $e$ as unit, $L$ a closed subgroup of $H$ and $(u_n)_n$ a sequence of $L$-biinvariant continuous functions of positive type on $H$.
\item[ {\rm (a)}] If $$\sum_{n=1}^\infty u_n(e)<\infty,$$
  then the series $\sum_{n=1}^\infty u_n$ converges uniformly on $H$
  and its sum is a $L$-biinvariant continuous function of positive type.
\item[ {\rm (b)}] Furthermore if, for $n\geq 1$,
  $$\sum_{k=1}^n u_k\ll\pp,$$ where $\pp$ is a $L$-biinvariant continuous function of
  positive type, then $$\sum_{n=1}^\infty u_n\ll\pp.$$
\item[ {\rm (c)}] If $v_n$ is another sequence such that $v_n\ll u_n$,
  then $$\sum_{n=1}^\infty v_n\ll\sum_{n=1}^\infty u_n.$$

\end{lem}

\begin{pro}\label{Pro8}For every subprojective system $\pp=\{\pp^{(k)}\}_k$ in $\Q$, there exists a projective system $\Phi=\{\Phi^{(k)}\}_k$ and functions $\psi^{(k)}\in\P^\bi(G(k))$ such that, for every $k$, \bq\label{eqthree}\pp^{(k)}=\Phi^{(k)}+\sum_{j=0}^\infty \res_k^{k+j}(\psi^{(k+j)}).\eq The functions $\Phi^{(k)}$ and $\psi^{(k)}$ are unique.
\begin{paragraph}{Proof.}\rm{Let $\pp\in\Q$. Put, for every $k\geq 1$, \bq\label{eqfour}\psi^{(k)}=\pp^{(k)}-\res_k^{k+1}(\pp^{(k+1)}).\eq By the definition of $\Q$, for every $k\geq 1$,
    $\psi^{(k)}$ is a function of positive type on $G(k)$. By iteration, equality (\ref{eqfour}) gives, for every $k\geq 1$, $$\pp^{(k)}=\psi^{(k)}+\res_k^{k+1}(\psi^{(k+1)})+\dots+\res_k^{k+n-1}(\psi^{(k+n-1)})+\res_k^{k+n}(\pp^{(k+n)}).$$ 
Put $\Psi^{(k,n)}=\sum_{j=0}^{n-1}\res_k^{k+j}(\psi^{(k+j)})$, then for every $k\geq 1$, $$\pp^{(k)}=\Psi^{(k,n)}+\res_k^{k+n}(\pp^{(k+n)}).$$ It follows that, for every $n\geq 
1$, $\Psi^{(k,n)} \ll \pp^{(k)}$. This implies, by (b) of Lemma \ref{Lem6},
that the sequence $\{\Psi^{(k,n)}\}_n$ converges uniformly on $G(k)$ to
$\Psi^{(k)}\in\P^\bi(G(k))$, where
$\Psi^{(k)}=\sum_{j=0}^\infty\res_k^{k+j}(\psi^{(k+j)}).$ Hence the
sequence $\res_k^{k+n}(\pp^{(k+n)})$ converges uniformly on $G(k)$. Let us denote by $\Phi^{(k)}$ its limit. Since $\res_k^{k+1}$
is continuous in the topology of uniform convergence on $G(k)$, \ba\Phi^{(k)}=\lim_{n\rightarrow+\infty}\res_k^{k+n}(\pp^{(k+n)})&=&\lim_{n\rightarrow+\infty}\res_k^{k+1+n}(\pp^{(k+1+n)})\\&=&\lim_{n\rightarrow+\infty}(\res_k^{k+1}\circ\res_{k+1}^{k+1+n})(\pp^{(k+1+n)})\\&=&\res_k^{k+1}\big(\lim_{n\rightarrow+\infty}\res_{k+1}^{k+1+n}(\pp^{(k+1+n)})\big)\\&=&\res_k^{k+1}(\Phi^{(k+1)}).\ea 

Then $\{\Phi^{(k)}\}_{k\geq 1}$ is a projective system. In order to prove the uniqueness, let us suppose that, for every $k\geq 1$, $\pp^{(k)}$ is given by another decomposition
$$\pp^{(k)}=\Phi_1^{(k)}+\sum_{j=0}^\infty\res_k^{k+j}(\psi_1^{(k+j)}),$$ then
\ba\psi^{(k)}&=&\pp^{(k)}-\res_k^{k+1}(\pp^{(k+1)})\\&=&\Phi_1^{(k)}+\sum_{j=0}^\infty\res_k^{k+j}(\psi_1^{(k+j)})\\&-&\res_k^{k+1}\bigg(\Phi_1^{(k+1)}+\sum_{j=0}^\infty\res_{k+1}^{k+1+j}(\psi_1^{(k+1+j)})\bigg)\\
&=&\sum_{j=0}^\infty\res_k^{k+j}(\psi_1^{(k+j)})-\sum_{j=1}^\infty\res_k^{k+j}(\psi_1^{(k+j)})=\psi_1^{(k)}.\hspace{1cm}\Box \ea }
\end{paragraph}
\end{pro} 

\begin{cor}\label{Cor2}Let $\pp_1=\{\pp_1^{(n)}\}_n$ and $\pp_2=\{\pp_2^{(n)}\}_n$ be two subprojective systems of $\Q$ such that $\pp_1\lll\pp_2$, in the sense that, for every $n$, $\pp_1^{(n)}\ll\pp_2^{(n)}$. Then, for every $n$, $\Phi_1^{(n)}\ll\Phi_2^{(n)}$ and $\psi_1^{(n)}\ll\psi_2^{(n)}$.
\begin{paragraph}{Proof.}\rm{We may write $$\pp_2=\pp_1+\pp_0, \ \rm{with} \ \pp_0\in\Q.$$ By the uniqueness of the decomposition given by formula (\ref{eqthree}), $$\Phi_2=\Phi_1+\Phi_0,$$ and for every $n$, $$\psi_2^{(n)}=\psi_1^{(n)}+\psi_0^{(n)}.$$ Since $\Phi_0^{(n)}$ and $\psi_0^{(n)}$ are in $\P^\bi(G(n))$, we can deduce that, for every $n$, $\Phi_1^{(n)}\ll\Phi_2^{(n)}$ and $\psi_1^{(n)}\ll\psi_2^{(n)}$.$\hspace{1cm} \Box$}
\end{paragraph}
\end{cor}

 By Corollary \ref{Cor1}, for every $n\geq 1$, $\P^\bi(G(n))$ is a lattice. Moreover, by Theorem \ref{Theo5}, $\P^\bi(G)$ is a lattice. Using the previous decomposition, we prove the following proposition, 

\begin{pro}\label{Pro9}The cone $\Q$ is a lattice.
\begin{paragraph}{Proof.}\rm{Let $\pp_1=\{\pp_1^{(n)}\}_n$, $\pp_2=\{\pp_2^{(n)}\}_n$ be two subprojective systems of $\Q$. By
    Proposition \ref{Pro8}, $$\pp_1^{(n)}=\Phi_1^{(n)}+\sum_{j=0}^\infty
    \res_n^{n+j}(\psi_1^{(n+j)}),$$
    $$\pp_2^{(n)}=\Phi_2^{(n)}+\sum_{j=0}^\infty
    \res_n^{n+j}(\psi_2^{(n+j)}).$$ Put, for every $n$, $\Phi_{Min}^{(n)}=\Phi_1^{(n)}\wedge\Phi_2^{(n)}$ and $\psi_{Min}^{(n)}=\psi_1^{(n)}\wedge\psi_2^{(n)}.$ Let
    $\pp=\{\pp^{(n)}\}_n\in\Q$. If $\pp\lll\pp_1$ and $\pp\lll\pp_2$, then by Corollary \ref{Cor2},
    for every $n$, $\Phi^{(n)}\ll\Phi_1^{(n)}$,
    $\Phi^{(n)}\ll\Phi_2^{(n)}$, and thus
    $\Phi^{(n)}\ll\Phi_{Min}^{(n)}$. Also, for every $n$,
    $\psi^{(n)}\ll\psi_1^{(n)}$, $\psi^{(n)}\ll\psi_2^{(n)}$, which
    implies that $\psi^{(n)}\ll\psi_{Min}^{(n)}$. Since, for every
    $n$, $\psi_{Min}^{(n)}\ll\psi_1^{(n)}$, then by (c) of Lemma \ref{Lem6},
    $\sum_{j=0}^\infty\res_n^{n+j}(\psi_{Min}^{(n+j)})$ converges in
    $\P^\bi(G(n))$ uniformly on $G(n)$. We put then, for
    every $n$, $$\pp_{Min}^{(n)}=\Phi_{Min}^{(n)}+\sum_{j=0}^\infty
    \res_n^{n+j}(\psi_{Min}^{(n+j)}).$$ We get, for every $n$,
    $\pp^{(n)}\ll\pp_{Min}^{(n)}$, and so $(\pp_1, \pp_2)$ has a greatest lower bound $\pp_{Min}=\{\pp_{Min}^{(n)}\}_n$. Now, put for every
    $n$, $\Phi_{Max}^{(n)}=\Phi_1^{(n)}\vee\Phi_2^{(n)},$ and $\psi_{Max}^{(n)}=\psi_1^{(n)}\vee\psi_2^{(n)}.$ Since, for
    every $n$, $\psi_{Max}^{(n)}\ll\psi_1^{(n)}+\psi_2^{(n)}$, then by (c)
    of Lemma \ref{Lem6}, we can put, for every $n$, $$\pp_{Max}^{(n)}=\Phi_{Max}^{(n)}+\sum_{j=0}^\infty \res_n^{n+j}(\psi_{Max}^{(n+j)}).$$ Thus, $(\pp_1, \pp_2)$ has a least upper bound $\pp_{Max}=\{\pp_{Max}^{(n)}\}_n$. As a consequence, $\Q$ is a lattice.$\hspace{1cm} \Box$}
\end{paragraph}
\end{pro}

Next, we will determine the set of extremal points of $\Q_{\leq1}$. We need to define, for $n\geq1$, the following subset :
$$\P^n=\{\pp\in 
\prod_{i=1}^\infty\P_{\leq1}^\bi(G(i))\ |\ \pp^{(i)}=\res_i^{n}\big(\pp^{(n)}\big),\ {\rm for}\ 1\leq i\leq n\ $$ $$ {\rm and}\ \pp^{(i)}=0,\ {\rm for}\ i\geq n+1\},$$ where, for every $i=1,\dots,n-1$, 
$$\res_i^{n}=\res_i^{i+1}\circ\res_{i+1}^{i+2}\circ\dots\circ\res_{n-1}^{n}.$$
The set $\P^n$, with finite $n$, consists of projective systems of finite order $n$ obtained via the following linear isomorphism : $$\iota :\P_{\leq1}^\bi(G(n))\rightarrow\P^n$$ 
$$\pp^{(n)}\longmapsto(\res_1^n(\pp^{(n)}),\res_2^n(\pp^{(n)}),\dots,\res_{n-1}^n(\pp^{(n)}),\pp^{(n)},0,\dots).$$

Since $\res_n^{n+1}(\P_{\leq1}^\bi(G(n+1)))\subset\P_{\leq1}^\bi(G(n))$, the set $\P_{\leq1}^\bi(G)$ can be identified with the projective limit of the family $\{\P_{\leq1}^\bi(G(n))\}_{n\geq 1}$ and an element $\pp\in\P_{\leq1}^\bi(G)$ determines a projective system $\{\pp^{(n)}\}$ with $\pp^{(n)}=\res_n(\pp)$. The same holds for an element $\omega$ of the set $\E_\infty$ of non zero extremal points of $\P_{\leq1}^\bi(G)$, i.e. $\E_\infty = {\rm Ext(}\P_{1}^\bi(G){\rm )}$. 

Let $\E_n$ denote the set of non zero extremal points of $\P^n$. An element $\pp\in\E_n$ is the image by the isomorphism $\iota$ of an element $\pp^{(n)}\in{\rm Ext(}\P_1^\bi(G(n){\rm )}$.

\begin{teo}The set of extremal points of $\Q_{\leq1}$ consists of  two types of elements : $$type \ \infty \ :\ \E_\infty,\  and \ type
  \ n\ :\ \E_n,$$ and we have \bq\label{eqfive} {\rm
    Ext(}\Q_{\leq1}{\rm
    )}={\rm \{}0{\rm\}}\cup\E_\infty\cup\big(\bigcup_{n=1}^\infty\E_n\big).\eq The sets
  $\E_\infty$, $\E_n$ $(n\geq 1)$ are disjoint.
\begin{paragraph}{Proof.}\rm{ (a) Let us prove that
    $\pp\in\E_n$ is extremal. Suppose that $\pp=\pp_1+\pp_2$, $\pp_1,
    \pp_2\in\Q_{\leq1}$. Then, for every $n$,
    $$\pp^{(n)}=\pp_1^{(n)}+\pp_2^{(n)}.$$ So,
    $\pp_1^{(n)}=\lambda_1\pp^{(n)}$,
    $\pp_2^{(n)}=\lambda_2\pp^{(n)}$. On the other hand, \ba\pp^{(n-1)}=\res_{n-1}^n\pp^{(n)}&=&\pp_1^{(n-1)}+\pp_2^{(n-1)}\\&\gg&\lambda_1\res_{n-1}^n\pp^{(n)}+\lambda_2\res_{n-1}^n\pp^{(n)}=\res_{n-1}^n\pp^{(n)}.\ea Therefore $$\pp_1^{(n-1)}=\lambda_1\res_{n-1}^n\pp^{(n)} ,\ \pp_2^{(n-1)}=\lambda_2\res_{n-1}^n\pp^{(n)},$$ and hence $$\pp_1=\lambda_1\pp , \ \pp_2=\lambda_2\pp.$$ 

(b) Let us prove that $\pp\in\E_\infty$ is extremal. Suppose that $\pp=\pp_1+\pp_2$, $\pp_1, \pp_2\in\Q_{\leq1}$. Since $\pp$ is a projective system, for every $n$, $\psi^{(n)}=0$. Thus, $\psi_1^{(n)}=0$, $\psi_2^{(n)}=0$, and hence $\pp_1, \pp_2\in\P_{1}^\bi(G)$. Therefore $$\pp_1=\lambda_1\pp ,\ \pp_2=\lambda_2\pp.$$

(c) Let $\pp$ be a non zero extremal point of $\Q_{\leq1}$, we can write $$\pp^{(n)}=\Phi^{(n)}+\sum_{j=0}^\infty \res_n^{n+j}(\psi^{(n+j)}),$$ it's a decomposition into two elements of $\Q_{\leq1}$:

\underline{First case} : $\psi^{(n)}=0$, for every $n$, and so $\pp\in\E_\infty$.

\underline{Second case} : $\Phi^{(n)}=0$, for every $n$, and hence $$\pp=\Psi_1+\Psi_2+\dots,$$ where \ba\Psi_n^{(j)}&=&\res_j^n(\psi^{(n)}) \ \ \ \ \ \ {\rm if} \ \ j\leq n,\\
               &=& \ \ \ \ \ \ \ 0 \ \ \ \ \ \ \ \ \ \ \ \ \ {\rm if}
               \ \ j>n.\ea As a result, there exists $n_0$ such that
               $\pp=\Psi_{n_0}$, with $\psi^{(n_0)}\in$
               Ext($\P_1^\bi(G(n_0))$). We can then conclude that $\pp\in\E_{n_0}$.$\hspace{1cm} \Box$ }
\end{paragraph}
\end{teo}

Assuming all $G(n)$ separable, we can now state a Bochner type theorem for the corresponding Olshanski spherical pairs.

\begin{teo}Let $(G,K)$ be an Olshanski spherical pair defined as inductive limit of an increasing sequence of Gelfand pairs $(G(n),K(n))_{n}$, with the assumption that all $G(n)$ are separable. Then, for every function $\pp\in\P^\bi(G)$, there exists, on the Borel set $\Omega=$\rm{
    Ext(}$\it\P_{\rm{1}}^\bi(G)$),\it \ a unique bounded and positive measure $\m$ such that  $$\pp(g)=\int_\Omega\omega(g)\m(d\omega).$$

\begin{paragraph}{Proof.}\rm{ The set $\Q_{\leq1}$ being
    convex, compact and metrisable in $\Q$, it satisfies the
    hypothesis of Choquet's theorem. Hence Ext($\Q_{\leq1}$) is a
    Borel set and every element of $\Q_{\leq 1}$ can be represented via a probability measure $\n$ on ${\rm Ext(}\Q_{\leq1}{\rm)}$ such that, for every continuous linear form $L$ on $\Q$, \bq\label{eqsix}
    L(q)=\int_{\rm{ Ext(}\Q_{\leq1}\rm{)}}L(p)\n(dp).\eq Moreover, as $\Q$ is a lattice (Proposition \ref{Pro9}), by (iii) of Choquet's theorem, the measure $\n$ is unique. Furthermore, we can deduce from formula (\ref{eqfive}) that $$\Omega={\rm Ext(}
    \Q_{\leq1}{\rm)}\setminus\big(\bigcup_{n=1}^\infty\E_n\cup {\rm\{}0{\rm\}}\big).$$ Hence $\Omega$ is a Borel set. 

Let $\pp$ be an element of $\Pin$. We know that $\pp$ determines a
sequence $\{\pp^{(n)}\}_{n\geq1}$ where $\pp^{(n)}=\res_n(\pp)$. Let us take, for $L$ in (\ref{eqsix}), the linear form $$\pp^{(n)}\mapsto (\pp^{(n)},f)=\int_{G(n)}\pp^{(n)}(h)f(h)\alpha_n(dh),$$ where $f\in L^1(G(n))$ and $\alpha_n$ is the left invariant Haar measure of $G(n)$. By considering, for every $n$, the approximation $(f_k)$ : $f_k\in L^1(G(n))$, $f_k\geq0$,$$\int_{G(n)}f_k(h)\alpha_n(dh)=1,$$ and for every continuous bounded function $\psi$ : $$\lim_{k\rightarrow\infty}\int_{G(n)}\psi(h)f_k(h)\alpha_n(dh)=\psi(g),$$ we get that, for every $n\geq 1$,
 $$\pp^{(n)}(g)=\int_{\Omega}\omega(g)\
 \n^{(\infty)}(d\omega)+\sum_{k=n}^\infty\int_{\E_n}\omega(g)\
 \n^{(k)}(d\omega),$$ where $\n^{(\infty)}$ (respectively $\{\n^{(k)}\}_{k\geq n}$), are the restrictions of $\n$ to $\Omega$ (respectively $\{\E_k\}_{k\geq n}$). Therefore we obtain, for $g\in G(n)$, $$\pp^{(n)}(g)-\pp^{(n+1)}(g)=\int_{\E_n}\omega(g)\
 \n^{(n)}(d\omega).$$ Since $\{\pp^{(n)}\}_{n\geq1}$ is a projective system, for every $g\in G(n)$ and every $n\geq 1$, $$\int_{\E_n}\omega(g)\
 \n^{(n)}(d\omega)=0.$$ As $\omega(e)=1$ we get, for every $n\geq 1$, $$\n^{(n)}(\E_n)=0.$$ Hence $\n$ is concentrated on $\E_\infty=\Omega$. It follows that every element $\pp\in\Pin$ has the following integral representation : $$\pp(g)=\int_\Omega\omega(g)\n^{(\infty)}(d\omega), \ {\rm(}g\in G{\rm)}.$$ Finally, every element $\pp\in\P^\bi(G)$ can be uniquely written as $\pp(g)=\lambda \pp_0(g)$ with $\pp_0\in\Pin$ and $\lambda=\pp(e)\geq 0$. Hence $\pp$ is represented via a measure $\m$ equal to $\lambda \n_0^{(\infty)}$, where $\n_0^{(\infty)}$ verifies $$\pp_0(g)=\int_\Omega\omega(g)\n_0^{(\infty)}(d\omega).\hspace{1cm}\Box$$ }
\end{paragraph}
\end{teo} 

\section{Remarks and open questions} (1) We do not know a topology
making $\Pin$ compact and enabling in consequence a direct application
of Choquet's theorem without using $\Q$. T. Hirai and E. Hirai had
studied this problem in \cite{Hir}.

(2) Given a generalized Gelfand pair, i.e. an Olshanski spherical
pair, one problem is to find the set of extremal points $\Omega$. This
is known in several cases. Another problem is, given $\pp\in\P^\bi(G)$, to find the representing measure $\m$.
\begin{paragraph}{Acknowledgements.}{\rm I would like to express my gratitude to my PhD advisors Professor Jacques Faraut and Professor Angela Pasquale for many extensive discussions and numerous suggestions during all stages of the preparation of this paper.}
\end{paragraph}

\end{document}